\documentclass[11pt,a4paper]{amsart}
\usepackage{amsmath}
\usepackage{amsfonts}
\usepackage{amssymb}
\newtheorem{theorem}{Theorem}[section]
\newtheorem{lemma}{Lemma}[section]
\newtheorem{proposition}{Proposition}[section]

\theoremstyle{definition}
\newtheorem{definition}{Definition}[section]

\theoremstyle{remark}
\newtheorem{remark}{Remark}[section]

\numberwithin{equation}{section}

\begin{document}

\title[Inv. Measure and the Euler Char. of Projective
Manifolds]{Invariant Measure and the Euler Characteristic of
Projectively Flat Manifolds}


\author{Kyeonghee Jo}
\address{School of Mathematical Science, Seoul National
University, 151-742 Seoul, KOREA} \email{khjo@math.snu.ac.kr}

\author{Hyuk Kim}
\email{hyukkim@math.snu.ac.kr}

\begin{abstract}
In this paper, we show that the Euler characteristic of an even dimensional closed
projectively flat manifold is equal to the total measure which is
induced from a probability Borel measure on $\mathbb RP^n$
invariant under the holonomy action, and then discuss its
consequences and applications. As an application, we show that the
Chern's conjecture is true for a closed affinely flat manifold
whose holonomy group action permits an invariant probability Borel
measure on $\mathbb RP^n$; that is, such a closed affinly flat
manifold has a vanishing Euler characteristic.
\end{abstract}
\thanks{*The author gratefully acknowledge the partial support of
BK 21 projection in 2002. **Research of this author was supported by
SNU-Daewoo and by grant No.R01-1999-00002 from the KOSEF}
\maketitle

\renewcommand{\thefootnote}{\alph{footnote}}
\setcounter{footnote}{0}

\section{Introduction}
\label{intro}
   In this paper, we will show how the Euler characteristic of a projectively
    flat manifold $M$ can be viewed as a total measure of $M$ where the measure
    is induced from a probability measure on $\mathbb RP^n$ invariant under the
    holonomy action, and then we will also discuss its consequences and
    applications.

   A projectively flat manifold $M$ is a manifold which is locally modelled
    on the projective space with its natural
    projective geometry, i.e, $M$ admits a cover of coordinate charts into
    the projective space $\mathbb RP^n$ whose coordinate transitions are
    projective transformations. By an analytic continuation of coordinate
    maps from its universal covering $\tilde M$, we obtain a developing map
    from $\tilde M$ into $\mathbb RP^n$ and this map is rigid in the sense
    that it is determined only by a local data. Therefore the deck
    transformation action on $\tilde M$ induces the holonomy action via
    the developing map by the rigidity.(See for example \cite{4,14,15} for more
    details.) Suppose there is a probability Borel measure $\lambda $ on $\mathbb RP^n$
    which is invariant under this holonomy action. Then we will first show
    that a Borel measure $\mu $ on $M$ is induced from $\lambda $ by the invariance
    property of $\lambda $, and then show the following Main Theorem.

\begin{theorem}[The Main Theorem]
 Let $M$ be an even dimensional closed projectively flat manifold and
 $ \lambda $ be a holonomy invariant finitely
 additive probability Borel measure on ${\mathbb R} P^n  $ .

 Then
$$ \chi(M) = \mu(M),$$
 where $\mu$ is the Borel measure on $M$ induced from $\lambda$.
\end{theorem}

  This result and its consequences are significantly refined and evolved
   versions of our earlier results \cite{8,9} in its perspective.
  Our investigations have been motivated from the effort to resolve the
   Chern's conjecture (or also known as Sullivan's conjecture):
  ``A closed affinely flat manifold has vanishing Euler
  characteristic.''

  An $(X,G)$-manifold is a manifold which is locally modelled on $X$ with
   the geometry determined by the Lie group $G$ acting on $X$ analytically.
  For example, projectively flat manifold is a special case of $(X,G)$-manifold
   with $X=\mathbb RP^n$ and $G$=PGL$(n+1,\mathbb R )$ and so is an affinely flat
   manifold with $X=\mathbb E^n$, the standard Euclidean space and $G=$ Aff($n$),
   the group of affine transformations on $\mathbb E^n$.
  An affinely flat manifold also can be viewed as a projectively flat manifold whose
   holonomy preserves the set of points at infinity, $\mathbb RP^{n-1}_{\infty }$,
   by identifying $\mathbb E^n$
   with the affine space given by $x_{n+1}=1$ in $\mathbb R^{n+1}$ so that
   $\mathbb RP^n$ becomes a compactification of $\mathbb E^n$.
  Similarly all the Riemannian and pseudo-Riemannian space forms can be
   considered as a subclass of projectively flat manifolds, and also a subclass
   of affinely flat manifolds if they are flat.

  The Euler characteristic of flat Riemannian or pseudo-Riemannian manifolds
   vanishes by Gauss-Bonnet-Chern theorem and
   it is natural to ask the same for
   affinely flat manifolds more generally and this is the content of Chern's
   conjecture.
  If a compact affinely flat manifold $M$ is complete, then the conjecture is true by
   the work of Kostant and Sullivan \cite{12},
    but note that the compactness does
   not necessarily imply the completeness in contrast to the Riemannian case.
  There has been various partial answers in different directions but the
   conjecture is not completely resolved yet.
  As one of the corollaries of the Main Theorem, we show that the conjecture is
   true if the holonomy group of affinely flat manifold has an invariant
   probability measure generalizing the earlier result for amenable case
   as well as for radiant case in a unified way.

  In Sect.~\ref{sec:1}, we will define a pull-back measure $f^* \lambda$ for a
   given local homeomorphism $f$ from a manifold $M$ to another manifold
   $N$ having a measure $\lambda$. For $(X,G)$-manifold
  $M$ and the corresponding developing map
    $D\, :\, \tilde M\to X$, $D^*\lambda $ is well-defined on $\tilde M$
    whenever $X$ has a measure $\lambda $.
   If $\lambda $ is invariant under its holonomy action, we can also define
    a measure $\mu $ on $M$ naturally induced from $D^*\lambda $ by covering
    projection $p\, :\, \tilde M\to M $.
   This is proved in Sect.~\ref{sec:2}.
   In Sect.~\ref{sec:3}, we will prove the Main Theorem.
   In Sect.~\ref{sec:4}, we will discuss the consequences and applications of the
   Main Theorem including the relation between the Euler
    characteristic and the developing maps for projectively flat manifolds
    and forementioned results for affinely flat manifolds.

\section{Pull-back measure}
\label{sec:1}
   Let $(N,\Omega,\lambda)$ denote a (finitely additive, resp.) measure space such
    that the $\sigma $-algebra (algebra, resp.) $\Omega $
    contains all the open sets of $N$,
    that is, $\lambda $ is a Borel measure if it is countably additive.
   In this paper we will also call such a measure finitely additive
    Borel measure  when $\lambda $ is only finitely additive.
   Let $M$ and $N$ be manifolds and
    $f\, :\, M\to N$ be a local homeomorphism.
   Then we can define a pull-back (finitely additive, resp.) Borel measure
    $f^* \lambda $ on $M$ whenever such a (finitely additive, resp.)  Borel measure
    $\lambda $ is given on $N$.

   An open covering $B=\{B_i\,|\, \overline{B_i}$ compact$\,
   \}_{i\in \mathbb N}$ of $M$ will be called a
    covering {\it adapted to local homeomorphism} $f$ if $f|_{B_i}$ is a
    homeomorphism for each $i$.

\begin{definition}\label{def:1.1}
 Let $B=\{B_i\}$ be an adapted covering. Define
\begin{equation*}
\Omega _B=\{\,A\subset M\, :\, f(A\cap B_i)\mbox{ is }\lambda
\mbox{-measurable for all }i\, \}.
\end{equation*}
\end{definition}

\begin{lemma}\label{lem:1.2}
 Let $B=\{B_i\}$ be an adapted covering.
 Then $\Omega _B$ is a $\sigma $-algebra.
 ($\Omega _B$ is an algebra if $ \lambda $ is finitely additive.)
\end{lemma}

\begin{proof}
   Obviously $\varnothing $ and $M$ are contained in $\Omega_B $.
   For each subset $A$ of $M$, we have the following equality:
\begin{equation*}
    f(A^c\cap B_i)=f|_{B_i}(B_i\setminus (A\cap B_i))
                   =f(B_i)\setminus f(A\cap B_i) \mbox{ for all }i.
\end{equation*}
    So the measurability of $f(A\cap B_i)$ implies the measurability of
     $f(A^c\cap B_i)$ for all $i$.
   Therefore $A^c\in \Omega_B$ whenever $A\in \Omega_B$.

   Let $\{S_j\}_{j\in I}\subset \Omega_B$ where $I$ is finite if $\lambda $
    is finitely additive and $I$ is equal to $\mathbb N$ if $\lambda $ is
    countably additive.
   Then $f\bigl ((\cup_{j\in I}S_j)\cap B_i\bigr )
         =f\bigl (\cup_{j\in I}(S_j\cap B_i)\bigr )
         =\cup_{j\in I}f(S_j\cap B_i)$ for all $i$.
   Therefore $f\bigl ((\cup_{j\in I}S_j)\cap B_i\bigr )$ is measurable for
    all $i$ and thus $\cup_{j\in I}S_j $ is contained in $\Omega_B$.

\end{proof}

\begin{definition}\label{def:1.3}
 Let $B=\{B_i\} $ be an adapted covering and $B_k'=B_k \setminus \cup_{i=1}^{k-1}B_i$.
 Define a set function
 $\lambda_B\, :\, \Omega_B \to \{x\in \mathbb R | x\geq 0\} \cup \{
 + \infty \} $ as follows.
$$ \lambda_B(A)=\sum_{k=1}^{\infty }\lambda (f(A\cap B_k'))$$
\end{definition}
$\lambda _B $ is  well-defined by the following reason. Firstly we see the following equalities
from the fact that $f$ is a homeomorphism on $B_j$ for each $j$:
\begin{equation*}
\begin{split}
f(B_k'\cap B_j)
 &=f|_{B_j}\bigl ((B_k \setminus \cup_{i=1}^{k-1}B_i)\cap B_j \bigr )\\
 &=f|_{B_j}\bigl ((B_k \cap B_j)\setminus (\cup_{i=1}^{k-1}B_i)\cap B_j \bigr )\\
 &=f(B_k \cap B_j)\setminus f((\cup_{i=1}^{k-1}B_i) \cap B_j )\mbox{ for each }j.
\end{split}
\end{equation*}
  Since $(\cup_{i=1}^{k-1}B_i )\cap B_j $ is open and $f|_{B_j}$ is a homeomorphism for each $j$,
  $f((\cup_{i=1}^{k-1}B_i )\cap B_j )$
  is also open and thus $\lambda $-measurable. This implies that $f(B_k'\cap B_j)$ is $\lambda $-measurable
  for all $j$ and so $B_k'\in \Omega_B$.
 So $A \cap B_k'\in \Omega_B$ for all $A\in \Omega_B $ and hence
 $f(A\cap B_k')=f(A\cap B_k'\cap B_k)$ is $\lambda $-measurable.
 Therefore $\lambda _B $ is  well-defined.

\begin{lemma}\label{lem:1.4}
 $(M,\Omega_B,\lambda_B)$ is a (finitely additive) Borel measure space.
\end{lemma}
\begin{proof}
   Let $ A$ be an open subset of $M$.
   Notice that $\Omega $ contains all the open sets of $N$ and hence
    $f(A \cap B_i) \in \Omega$.
   Therefore all the open sets of $M$ are contained in $ \Omega_B.$
   Now, it suffices to show the countable (finite) additivity of $\lambda $.
   Let $\{S_i\}_{i\in I}\subset \Omega_B$ be a disjoint collection,
    where $I$ is finite if $\lambda $ is finitely additive and $I$ is equal
    to $\mathbb N$ if $\lambda $ is countably additive.
   Then we have
\begin{equation*}
\begin{split}
\lambda _B(\cup_{i\in I}S_i) &=\sum_{k=1}^{\infty }\lambda \bigl
(f((\cup_{i \in I}S_i)\cap B_k')\bigr )\\
                             &=\sum_{k=1}^{\infty }\lambda \bigl (\cup_{i\in
                               I}f|_{B_k}(S_i\cap B_k')\bigr )\\
                             &=\sum_{k=1}^{\infty }\sum_{i\in I}\lambda
                               \bigl (f(S_i\cap B_k')\bigr )\\
                             &=\sum_{i\in I}\sum_{k=1}^{\infty }\lambda
                               \bigl (f(S_i\cap B_k')\bigr )\\
                             &=\sum_{i\in I}\lambda_B(S_i).
\end{split}
\end{equation*}
   Notice that the third equality holds since $f|_{B_k}$ is a homeomorphism.
\end{proof}

\begin{lemma}\label{lem:1.5}
 Let  $\{B_i\}$ and $\{C_i\}$ be adapted coverings and
 $\{\Omega_B, \lambda_B\}$ and $\{\Omega_C,\lambda_C\}$ be the corresponding
 measure systems on $M$.
 Then $\Omega _B=\Omega_C$ and $\lambda_B=\lambda_C$.
\end{lemma}
\begin{proof} (i) $\Omega_B=\Omega_C$:
   For each $C_k$ and for all $i$, $f(C_k\cap B_i)$ is open and thus
    $f(C_k\cap B_i)$ is $\lambda $-measurable.
   Hence $C_k\in \Omega_B$ for all $k$.
   Suppose $A\in \Omega_B$. To prove $A\in \Omega_C$, it suffices to show
    $f(A \cap C_i)$ is $\lambda $-measurable for all $i$.
\begin{equation*}
\begin{split}
    f(A\cap C_i) &=f\bigl ((A\cap C_i)\cap (\cup_{k=1}^{\infty }B_k)\bigr )\\
                 &=f\bigl (\cup_{k=1}^{\infty }(A\cap C_i\cap B_k)\bigr ) \\
                 &=\cup_{k=1}^{\infty }f(A\cap C_i\cap B_k)
\end{split}
\end{equation*}
   The last union is in fact a finite union
    since $\overline {C_i}$ is compact.
   Since $C_i\in \Omega_B $ for all $i$, we get $A\cap C_i\in \Omega _B$
    for all $i$
   and so $f(A\cap C_i\cap B_k)$ is $\lambda $-measurable for   all $i,\, k$.
   Therefore  $f(A\cap C_i)$ is measurable
    whether $\lambda $ is finitely additive or not.
   This shows $\Omega_B \subset \Omega_C $ and similarly
    $\Omega_C \subset \Omega_B $
\newline
(ii) $\lambda_B=\lambda_C$:
    For each $A\in \Omega_B=\Omega_C$, we have
\begin{equation*}
\begin{split}
    \lambda_B(A) &=\sum_{k=1}^{\infty }\lambda \bigl (f(A\cap B_k')\bigr )\\
                 &=\sum_{k=1}^{\infty }\lambda \bigl (f((A\cap B_k')\cap
                   (\cup_{j=1}^{\infty }C_j')\bigr )\\
                 &=\sum_{k=1}^{\infty }\lambda \bigl (\cup _{j=1}^{\infty }
                   f(A\cap B_k'\cap C_j')\bigr )\\
                 &=\sum_{k=1}^{\infty }\sum_{j=1}^{\infty }\lambda
                   \bigl (f(A\cap B_k'\cap C_j')\bigr )\\
                 &=\sum_{j=1}^{\infty }\bigl (\sum_{k=1}^{\infty }\lambda
                   (f(A\cap C_j'\cap B_k'))\bigr )\\
                 &=\sum_{j=1}^{\infty }\lambda \bigl (f(A\cap C_j')\bigr )\\
                 &=\lambda_C(A),
\end{split}
\end{equation*}
   where the fourth, fifth and the last  equalities hold
    even if $\lambda $ is finitely additive
    since  the summation $\sum_{j=1}^{\infty }$ and
    $\sum_{k=1}^{\infty }$ are in fact finite sums.
\end{proof}
   From Lemma~\ref{lem:1.4} and~\ref{lem:1.5}, we have the following theorem.

\begin{theorem}\label{thm:1.6}
 Let $f\,:\, M \to N$ be a local homeomorphism and
 $\lambda $ a (finitely additive, resp.) Borel
 measure on $N$.
 Then there exist a (finitely additive, resp.) measure $\mu $ on $M$
 such that
 \begin{enumerate}
  \item [\rm(i)] every open subset of $M$ is measurable,
  \item [\rm(ii)] for $\mu $-measurable  subset $A$ of $M$ (with $\overline A$
   compact, resp.)
   $\mu (A)= \lambda (f(A))$ if $f|_A$ is a homeomorphism.
 \end{enumerate}
 Furthermore, $\mu $ is unique in the following sense:
 If $\mu '$ is another measure satisfying the above property, then
 $\mu (A) =\mu '(A)$ for each open set $A$ (with $\overline A$ compact,
 resp.) of $M$.
\end{theorem}
\begin{proof}
   The existence of such a measure is proved by Lemma~\ref{lem:1.4},~\ref{lem:1.5} and
    the definition of $\mu $ respectively.
   So the only thing left to prove is the uniqueness.
   Let $A$ be a open subset of $M$ and $\{B_i \}$ be an adapted covering.
   Then $A$ and $B_i$ are $\mu'$-measurable and
        $ \mu (A) =\sum_{k} \mu (A \cap B_k')
                  =\sum_{k} \lambda(f (A \cap B_k'))
                  =\sum_{k} \mu' (A \cap B_k')
                  =\mu'(\cup_k (A \cap B_k'))
                  =\mu'(A). $
   Notice that  $\sum_{k}$ is in fact a finite sum if $\overline A$ is
    compact and thus the first and fourth equalities also hold in
    the finitely additive case.
   Also the second and third equalities hold by the property (ii).
\end{proof}

   Up to now we have defined a pull-back Borel measure on $M$ when  a local
    homeomorphism $ f\,:\, M \to N$ and a Borel measure $\lambda $
    on $N$ are given.
   Now, we will denote the pull-back measure by  $f^*\lambda $.
   Notice, for any measurable subset $A$ of $M$,
    $f^*\lambda(A) = \sum_i \lambda \bigl (f(A\cap U_i')\bigr )$
    where $\{U_i\}$ is any adapted covering of $M$ and
    $U_k'=U_k \setminus \cup_{i=1}^{k-1}U_i.$
   Denote the $\sigma $-algebra corresponding to $f^*\lambda $ by
    $f^*\Omega $.

\begin{remark}
   If $f$ is a homeomorphism, then $f^*\Omega =f^{-1}\Omega $ and
    $f^*\lambda (A)=\lambda (f(A))$ for all $A\in f^*\Omega $ when $\lambda $
    is countably additive.
   But this does not hold if $\lambda $ is finitely
    additive.
   In fact, id$^*\lambda \neq \lambda $ if $\lambda $ is finitely additive
    and $M$ is not compact. For example,
    if $M=\mathbb R$ and $\lambda $ is any finitely additive translation
    invariant probability measure of $\mathbb R$ (the amenability of $\mathbb R$
    ensures the existence of such a measure),
    then any bounded measurable subset of $\mathbb R$ has a measure $0$ and thus
    id$^* \lambda \equiv 0$ by the definition of id$^* \lambda $.
   This strange phenomenon arises since our measure id$^* \lambda$ is pulled
    back only locally and then is given as the sum of these local measures
    not reflecting the global nature  of the original measure $\lambda $.
\end{remark}
\begin{theorem}\label{thm:1.7}
 Suppose topological groups $G$ and $H$ act continuously on $M$ and $N$ respectively.
 Let $(\phi ,f)\,:\, (G,M) \to (H,N)$ be an equivariant pair where
 $\phi \,:\,G\to H$ is a homomorphism and $f\,:\,M\to N$ is a local
 homeomorphism, and $\lambda $ be an $H-$invariant Borel measure.
  Then $f^*\lambda$ is $G-$invariant.
\end{theorem}
\begin{proof}
   Let $A$ be a measurable subset of $M$ and $\{B_i\}$ be an
    adapted covering of $M$.
   Then for each fixed $g \in G$, $\{gB_i\}$ is also an adapted covering of $M$.
\begin{equation*}
\begin{split}
    f^*\lambda (gA)&=\sum \lambda \bigl (f(gA\cap (gB_i)')\bigr )\\
                   &=\sum \lambda \bigl (f(gA\cap gB_i')\bigr )\\
                   &=\sum \lambda \bigl (f(g(A\cap B_i'))\bigr ) \\
                   &=\sum \lambda \bigl (\phi (g)f(A\cap B_i')\bigr )\\
                   &=\sum \lambda \bigl (f(A\cap B_i')\bigr )\\
                   &=f^*\lambda (A)
\end{split}
\end{equation*}
   Notice  the fifth equality holds since $\lambda $ is $H-$invariant.
\end{proof}

\section{Holonomy invariant measure}\label{sec:2}
   Let $p\, :\, \widetilde M\to M$ be a regular covering map and $\lambda $
    be a Borel measure on $\widetilde M$.
   Assume $\lambda $ is invariant under the action of the deck transformation
    group.
   Then we will define a Borel measure $\mu $ on $M$ such that $p^*\mu =\lambda $.

   An open covering $B=\{B_i\}_{i\in \mathbb N}$ of $M$ will be called a
    covering {\it adapted to covering map $p$},
    if $B_i$ is evenly covered and $\overline {B_i}$ is compact for each $i$.

\begin{definition}\label{def:2.1}
 Let $B=\{B_i\}$ be an adapted covering. Define
\begin{equation*}
\begin{split}
\Omega_B =\{& A\subset M\, |\mbox{ for each }i,\bigl
(p|_{\widetilde {B_i}}\bigr )^{-1}(A\cap B_i)\mbox{ is }\lambda
\mbox{-measurable on }\widetilde M\\ &\mbox{for some lifting
}\widetilde {B_i}\, \}
\end{split}
\end{equation*}
\end{definition}
   Note that since $p$ is regular and $\lambda $ is invariant under the
    action of the deck transformation group, $A \in \Omega_B$ implies
    $\bigl (p|_{\widetilde {B_i}}\bigr )^{-1}(A\cap B_i)$ is
    $\lambda $-measurable for any lifting $\widetilde {B_i}$, that is
    $\Omega _B$ is well-defined independently of lifting.
   The following Lemmas~\ref{lem:2.2}, ~\ref{lem:2.4} and ~\ref{lem:2.5} are easily proved by  the same
    argument as in Sect.~\ref{sec:1} and thus we will omit the proofs.
\begin{lemma}\label{lem:2.2}
 Let $B=\{B_i\}$ be an adapted covering of $M$.
 Then $\Omega_B$ is a $\sigma $-algebra.
 ($\Omega_B$ is an algebra if $\lambda $ is finitely additive.)
\end{lemma}

\begin{definition}\label{def:2.3}
 Let $B=\{B_i\}$ be an adapted covering of $M$.
 Define a function
\begin{equation*}
\mu _B : \Omega _B \to \{x\in \mathbb R \, :\, x \geq 0\} \cup \{
 + \infty \} \mbox{ by }
\mu_B(A)=
 \sum_{k=1}^{\infty }
 \lambda \bigl (\bigl (p|_{\widetilde {B_k}}\bigr )^{-1}(A\cap B_k')\bigr ),
\end{equation*}
where $B_k'=B_k\setminus \cup_{i=1}^{k-1}B_i$.
\end{definition}
   It is easy to prove $B_k'\in \Omega_B$ by similar argument as in Sect.~\ref{sec:1}.
   Hence, $\mu_B$ is well defined independently of the choice of the lifting
    $\widetilde {B_i}$ since $\lambda $ is invariant under the action of the
    deck transformation group.
\begin{lemma}\label{lem:2.4}
 Let $B=\{B_i\}$ be an adapted covering of $M$.
 Then $(M,\Omega_B,\mu_B)$ is a (finitely additive) Borel measure space.
\end{lemma}

\begin{lemma}\label{lem:2.5}
 Let both $\{B_i\}$ and $\{C_i\}$ be  adapted
 coverings of $M$ and $\{\Omega_B,\mu_B\}$ and $\{\Omega_C,\mu_C\}$ be
 the corresponding measure systems on $M$.
 Then $\Omega_B=\Omega_C$ and $\mu_B =\mu_C$.
\end{lemma}

   From the above Lemmas, we have the following theorem.
\begin{theorem}\label{thm:2.6}
 Let $p\, :\, \widetilde M\to M$ be a regular covering and $\lambda $ be a
 (finitely additive, resp.) Borel measure on $\widetilde M$.
 Assume that $\lambda $ is invariant under the action of the deck
 transformation group.
 Then there exists a (finitely additive, resp.) Borel measure $\mu $ on $M$
 such that  $p^*\mu = \lambda $.
 Furthermore $\mu $ is unique in the following sense:
 if ${\mu }'$ is another (finitely additive, resp.) Borel measure
 such that  $p^*{{\mu }'}= \lambda $,
 then $\mu (A) = {\mu }^{\prime }(A)$ for each open subset $A$
 (with $\overline A$ compact, resp.) of $M$.
\end{theorem}
\begin{proof}
   The existence is obvious by Lemmas~\ref{lem:2.4} and ~\ref{lem:2.5}.
   Let $A$ be an open subset of $M$ and $\{B_i \}$ be an adapted covering.
   Choose a lifting $\widetilde {B_i}$ for each $i$.
   Then $\mu(A)=\sum_{k}\mu (A \cap B_k')
               = \sum_{k}p^* \mu' \bigl (\bigl
                 (p|_{\widetilde {B_k}}\bigr )^{-1}(A\cap B_k') \bigr )
               = \sum_{k}\mu'(A\cap B_k')= \mu' (A)$.
  Notice that the summation $\sum_k$ is in fact a finite sum if $\overline A$
   is compact and thus the first and fifth equalities also hold in
   the finitely additive case.
  Also the second and third equalities hold since $\overline {A \cap B_k'}$
   is compact.
\end{proof}

\begin{theorem}\label{thm:2.7}
 Let $M$ be a $(X,G)-$manifold,  $D\, :\,\widetilde M \to X$ be a developing
 map with holonomy group $H$, $p\, :\, \widetilde M \to M$ be a covering map
 and $\lambda $ be an $H-$invariant (finitely additive, resp.) Borel measure
  on $X$.
 Then there exists a (finitely additive, resp.) Borel measure $\mu $ on $M$
 such that $p^*\mu =D^* \lambda $.
\end{theorem}
\begin{proof}
   By Theorem~\ref{thm:1.7}, $D^* \lambda $ is invariant under the action of deck
    transformation since $\lambda $ is $H$-invariant.
   So there exists a Borel measure $\mu $ on $M$ such that
   $p^*\mu =D^* \lambda$  by Theorem~\ref{thm:2.6}.
\end{proof}

   The measure $\mu $ described in Theorem~\ref{thm:2.7} will be called the
    {\it {induced measure}}.
   Recall that a $(G,X)-$manifold is a smooth manifold which has a cover of
    coordinate charts by open subsets in $X$ whose coordinate transition
    functions are locally restrictions of the elements of $G$.
   To distinguish from topological chart, we will call a coordinate chart
    constituting a $(G,X)-$manifold as a {\it {geometric chart}}.

\begin{remark}
   Let $A$ be a subset of $M$ so that $A$ is contained in some evenly covered
    geometric chart.
   Assume that $\overline A $ is compact if $\lambda $ is finitely additive.
   Then $\mu(A)=\lambda(D\tilde A)$ where $D\tilde A$ is any developing image of a lifting $\tilde A$ of $A$.
\end{remark}

\section{Proof of the main theorem}\label{sec:3}
   We use a generalized {\it {Gauss-Bonnet formula}} in terms of angles of
    simplices in a triangulation to prove the theorem as in \cite{9}.
   The notion of angle is not well defined in general, but we can do define
    an angle using an holonomy invariant measure.

   Let (${\mathbb R}P^n$, $\lambda$) be a finitely additive probability
    Borel measure.
   We will use the same symbol $\lambda $ to denote the pull-back measure
    induced on ${\mathbb S}^n$ which is invariant under the antipodal map so that
    $\lambda  ({\mathbb S}^n)=2$.

   Let $s^n$ be a spherical simplex lying in the standard unit sphere
    ${\mathbb S}^n\subset {\mathbb R}^{n+1}$ so that each of its
    ($n$-$1$)-dimensional faces is a part of great hyperplanes
    $P_i,$ $ i=1,2,\cdots$,$n$+$1$, in general position.
   A great hyperplane is the intersection of ${\mathbb S}^n$ with an $n$-dimensional subspace
    of ${\mathbb R}^{n+1}$.
   Let $f_i$ be the characteristic function of {\it positive} side of $P_i$
    which, by definition, is the half of ${\mathbb S}^n$ bisected by $P_i$ that
    contains $s^n$.

   If $s^r < s^n$ is an $r$-dimensional face of $s^n$ given by
    $s^r = P_{i_1} \cap \dots \cap P_{i_{n-r}} \cap s^n$, the {\it angle} at
    $s^r$ in $s^n$, denoted by $\alpha (s^r,s^n)$, is defined as
\begin{equation*}
\alpha(s^r,s^n)
   = \frac{1}{2} \int_{S^n} f_{i_1} \cdots f_{i_{n-r}} d\lambda
\end{equation*}
   Clearly,
\begin{equation*}
\begin{split}
   \lambda (s^n) & = \int_{S^n} f_1 \cdots f_{n+1} \,d\lambda \\
   \lambda  ({\overline s}^n) & =  \int_{S^n} (1-f_1)
   \cdots (1-f_{n+1}) \,d\lambda
\end{split}
\end{equation*}
   where ${\overline s}^n$ is the antipodal image of $s^n$.
   Then
\begin{equation*}
\begin{split}
   \lambda (s^n) = \lambda ({\overline s}^n) & = \int_{S^n} (1-f_1) \cdots
   (1-f_{n+1}) \, d\lambda    \\
   & = \sum_{r=0}^n \sum_{s^r < s^n} (-1)^{n-r} 2\alpha (s^r,s^n) +
   (-1)^{n+1}\int_{S^n}  f_1 \dots f_{n+1} \, d\lambda
\end{split}
\end{equation*}
   and hence we get the {\sl Spherical Gauss-Bonnet formula} \cite{7}
\begin{equation*}
2\sum_{r=0}^n \sum_{s^r < s^n} (-1)^{n-r} \alpha (s^r,s^n)
   = (1+(-1)^n) \lambda (s^n).
\end{equation*}
or equivalently,
\begin{equation*}
2\sum_{r=0}^n \sum_{s^r < s^n} (-1)^r \alpha (s^r,s^n)
   = (1+(-1)^n) \lambda (s^n).
\end{equation*}
   Let $k(s^n) = \sum_{r=0}^n \sum_{s^r < s^n} (-1)^r \alpha (s^r,s^n)$.
   Then the above formula implies $k(s^n)$ equals zero when $n$ is odd and
    equals $\lambda (s^n)$ when $n$ is even.

   Let $M$ be a closed projectively flat manifold with a {\sl geometric
    triangulation} $K$ consisting of simplices {${\sigma}^n$} whose developing images are
    spherical simplices.
   Let $D:\tilde M \rightarrow {\mathbb R}P^n$ and $H$ be the corresponding
    developing map and holonomy group respectively.
   Let  $\lambda$ be an $H$-invariant finitely additive probability measure
    on ${\mathbb R}P^n$.
   Then for each face ${\sigma}^r < {\sigma} ^n$ the  {\it angle} at
    ${\sigma}^r$ in ${\sigma}^n$, denoted by
    $\tilde {\alpha }({\sigma }^r,{\sigma }^n)$, is defined as
\begin{equation*}
\tilde {\alpha}({\sigma}^r,{\sigma}^n) = \alpha(s_0^r,s_0^n),
\end{equation*}
  where $(s_0^r,s_0^n)$ is a developing image of a lifting of a pair
  $({\sigma}^r,{\sigma}^n)$ in $\tilde M$.

Note that
   $\tilde {\alpha}({\sigma}^r,{\sigma}^n)$ is well defined by the following
    reason.
   $\alpha(s_0^r,s_0^n)$ is actually the measure of some subset of
    $\mathbb RP^n$ and for any other choice of a lifting and its developing
    image $s_1$ there exists $h \in H$ so that
    $s_1^n = h(s_0^n)$ and $s_1^r = h(s_0^r)$.
   Therefore we get $\alpha (s_0^r,s_0^n) = \alpha (h(s_0^r),h(s_0^n))
    = \alpha (s_1^r,s_1^n)$ since $\lambda $ is $H$-invariant.
   But the angle depends on the developing map.

   From now on, we'll simply denote the angle $\tilde{\alpha}$ by $\alpha$
    by abusing notation

   Let
\begin{equation*}
\begin{split}
   S({\sigma}^i) & = \sum_{{\sigma}^i < {\sigma}^n}
                     \alpha ({\sigma}^i,{\sigma}^n),\\
   k({\sigma}^n) & = \sum_{r=0}^n (-1)^r \sum_{{\sigma}^r < {\sigma}^n}
   \alpha ({\sigma}^r,{\sigma}^n),
\end{split}
\end{equation*}
and  for a vertex $\nu$  in $K$, let
\begin{equation*}
d(\nu)  =
\sum_{r=0}^n \frac {(-1)^r}{r+1} \sum_{\nu \in {\sigma}^r}
  (1-S({\sigma}^r)).
\end{equation*}
   Then it is basically a rearragement of the angle terms to verify the
    following {\sl polyhedral Gauss-Bonnet formula} for the Euler
    characteristic $\chi(M)$:
\begin{equation*}
\sum_{\nu \in K} d(\nu) + \sum_{{\sigma}^n \in K} k({\sigma}^n) =
\chi(M).
\end{equation*}
   See \cite{8} for a proof and see \cite{10} for a motivation and geometric
    meanings of the terms in the formula.

   Let PGL$(n+1,\mathbb R)$, as usual, be the projective general linear group,
    i.e, GL$(n+1,\mathbb R)/ {\mathbb R^* I}$, where $\mathbb R^*$ is
    the set of nonzero elements of $\mathbb R$. If $V,W \subset \mathbb R ^{n+1}$ is
    a nonzero linear subspace, we denote by $\lbrack V\rbrack \subset
    \mathbb RP^n$ its image in $\mathbb RP^n$ and write $W\lneq V$ if $W$ is a
    proper subspace of $V$.

   Let $S=\lbrace V\lneq \mathbb R ^{n+1}\, |\, \lambda(\lbrack V \rbrack )>0 $
    and
  $ \lambda(\lbrack W \rbrack )=0$  for any $W \lneq V \rbrace$.
      Then $S$ is a disjoint union of $S_i$ with $i\in \lbrace
      0,1,\cdots, n-1\rbrace$, where
\begin{equation*}
S_i = \lbrace V\in S \,|\, \text{dim} V = i+1\rbrace.
\end{equation*}
 Let $S_{i,j} = \lbrace V\in S_i \,|\, \lambda([V])>
\frac 1j\rbrace . $ Then $|S_{i,j}|<j$, since $\lambda(\mathbb
RP^n)=1$ and $\lambda ([V \cap W])=0$ if $V$ and $W$ belong to
$S_{i,j}$ and $V \neq W$. Therefore $S_i =\cup_{j=1}^{\infty}S_{i,j}$
is countable and so is $S=\cup_{i=0}^{n-1}S_i$. Therefore we have
the following properties:
\begin{enumerate}
   \item [\rm (i)] $|S| $ is countable.
   \item [\rm (ii)] If $X\lneq \mathbb R^{n+1}$ and $X$ is transversal to
   each element of $S$,
   then $\lambda (\lbrack X \rbrack )=0$.
\end{enumerate}

   Therefore (i) implies that we can choose a geometric triangulation $K$
   on $M$ by a small perturbation
   such that every hyperplane in $\mathbb RP^n $ containing a developing
   image of some $(n-1)$-dimensional geometric simplex in $K$ is transversal
   to $S$ and so it has a measure zero by (ii)

   We now prove the Main Theorem:
    $S({\sigma}^i)=1$ for any geometric simplex ${\sigma}^i$ in $K$ by the above consideration and thus
    $d(\nu)=0$ for all vertex $\nu$ $\in$ $K$.
   Let $s_0^n$ be any developing image of ${\sigma}^n$.
   Then $k({\sigma}^n) = k(s_0^n)$  by definition of
    $\alpha({\sigma}^r,{\sigma}^n)$ and thus we get
    $k({\sigma}^n) = k(s_0^n) = \lambda(s_0^n)$.
   We may assume ${\sigma}^n$ is evenly covered and lies in some geometric
    chart, $\lambda(s_0^n)=\mu({\sigma}^n)$ by Remark in Sect.~\ref{sec:2}.
   Therefore $k({\sigma}^n) = \mu({\sigma}^n)$ for all
    $n$-simplex ${\sigma}^n$ $\in$ $K$.
   Now by the {\sl polyhedral Gauss-Bonnet Theorem},
\begin{equation*}
\chi(M) = \sum_{{\sigma}^n \in K} \mu({\sigma}^n)
\end{equation*}
   But we have chosen a triangulation so that
    the faces of $\sigma ^n$ have measure zero
    and hence
\begin{equation*}
\sum_{{\sigma}^n \in K} \mu({\sigma}^n) = \mu(M).
\end{equation*}
   This completes the proof.\qed

\section{Consequences and Applications}\label{sec:4}
   The right hand side $\mu (M)$ of the formula in the Main Theorem is
    supposed to depend on the holonomy invariant measure chosen and on
    the developing map, namely the projectively flat structure of $M$.
   But the theorem says that in fact it does not, and is always equal
    to the Euler charateristic of $M$, a topological invariant.
   Futhermore, there is no reason, a priori, that the total measure of
    $M$, $\mu (M)$ should be an integer.
   The topology, geometry and the measure related to $M$ are interlocked
    by the formula and we can expect interesting applications from these
    observations.
   We will see some of the immediate consequences and applications in
    this section.

   Let $M$ be a closed projectively flat manifold with amenable holonomy
    group $H$ and $m$ an invariant mean on $B(H)$, the space of all bounded functions
    on $H$, see \cite{5} for definitions of amenable group and invariant mean. We may assume that
    $m$ is right invariant since $H$ is a group, that is, $m(f_s)=m(f)$ for all $s \in H$, where $f_s$
    is a bounded function on $H$ given by $f_s(t)=f(ts)$.
   Then we can define an $H$-invariant finitely additive probability measure
    on ${\mathbb R}P^n$ as follows.
   Choose any probability measure ${\lambda}_0$ on ${\mathbb R}P^n$.
   Then for any ${\lambda}_0$-measurable subset $E$ of ${\mathbb R}P^n$ we can
    define a bounded function $f_E: H \rightarrow [0,1]$ by
\begin{equation}\label{eq:4.1}
 f_E(h) = \lambda_0(h(E))
\end{equation}
    for all $h \in H$.
   Now define a new measure $\lambda$ on ${\mathbb R}P^n$ by
\begin{equation*}
\lambda(E) = m(f_E)
\end{equation*}
    for all ${\lambda}_0$-measurable subset $E \subset {\mathbb R}P^n$.
   Then, by the property of invariant mean, $\lambda$ is a finitely additive
    $H-$invariant probability  measure on ${\mathbb R}P^n$.
More precisely,  $$\lambda(hE)
 =m(f_{hE})=m((f_E)_h)=m(f_E)=\lambda(E)$$
 since $$f_{hE}(h')=\lambda_0(h'(hE))=\lambda_0(h'h)E)=(f_E)(h'h)=(f_E)_h(h'),$$
 and the property
 $m(1)=1$ implies that $\lambda$ is a finitely additive probability  measure on ${\mathbb R}P^n$.

\begin{theorem}\label{thm:4.1}
 Let $M$ be an even dimensional closed projectively flat manifold with
 holonomy group $H$. Suppose there exist an $H$-invariant finitely additive probability Borel
 measure $\lambda $ on $\mathbb RP^n$. Then we have the following.
\begin{enumerate}
 \item [\rm (i)]   $\chi(M)$ is nonnegative.
 \item [\rm (ii)]  If the developing map is injective,
              then $\chi(M) = 0$.
 \item [\rm (iii)] If $\chi(M) =0$, then the developing map is not surjective.
 \item [\rm (iv)]  If the invariant measure $\lambda $ is countably additive,
              then $\chi(M)=0$ if and only if $\lambda(\Omega)=0$ and $\chi(M)>0$
              if and only if $\lambda(\Omega)=1$, where $\Omega$ is the developing image.
\end{enumerate}
\end{theorem}
\begin{proof}
   To prove this theorem, it suffices to prove (ii), (iii) and (iv)
    because (i) is the immediate consequence of the Main Theorem.

   (ii) Let $F$ be an open fundamental domain of $M$ in $\tilde {M}$.
   Let $D:\tilde M\rightarrow {\mathbb R}P^n$  be the corresponding developing
    map and $\phi :\pi_1(M)\rightarrow H$ be the holonomy representation.
   Suppose $\phi (\xi )=$ id for some $\xi \in \pi _1(M)$. Then $D(\xi x)=
    \phi (\xi) D(x) =D(x)$ for all $x\in \tilde M$. Since $D$ is injective,
    $\xi x=x $ for all $x\in \tilde M $, i.e, $\xi =$ id.
   Therefore $\phi $ becomes an isomorphism.
   Note that $H$ is non-trivial: If it were trivial, $M$ is simply connected
    and $D$ becomes a homeomorphism since $M$ is compact and $D$ is an
    injective local homeomorphism.
   But this is impossible since
    $\pi _1({\mathbb R}P^n)= \mathbb Z /2 $ for $n \geq 2$. Let $h\in H$ be a
    non-identity element.
   Then by the injectivity of $D$, $D(F)\cap h(D(F))=\varnothing $.
   Since $\lambda(D(F))=\lambda(h(D(F)))$ and
   $\lambda (D(F))+\lambda (h(D(F))) \leq 1$, the Euler characteristic,
   being a non-negative integer, has to be zero.

   (iii) Let $\mu $ be the induced measure. If $\chi (M)=0$,
   then the Main Theorem implies that $\mu (M)=0$ and thus for each $x\in D(
   \tilde M)$ there exists an open neighborhood $U_x$ of $x$
    such that $\lambda (U_x)=0$ by the definition of the induced measure.
   Therefore any compact subset $E$ in $D(\tilde M)$ has measure
    zero.
   Suppose $D$ is surjective. Then $\mathbb RP^n=D(\tilde M)$.
   But $\mathbb RP^n$ is compact and thus $\lambda (\mathbb RP^n)=0$.
   This contradicts that $\lambda$ is a probability measure.

   (iv) Let $\mu$ be the induced measure on $M$.
   $\chi (M)=0$ implies $\mu (M)=0$ by the Main Theorem
    and thus again for each $x \in \Omega$ there exists an open
    neighborhood $U_x$ of $x$ such that $\lambda (U_x)=0 $ by definition of
    the induced measure.
   Therefore any compact subset $E$ in $\Omega $
    has measure zero and thus $\lambda (\Omega ) =0$ by the countable
    additivity of $\lambda $.
   The converse is clear. Therefore $\chi (M)>0$ implies $\lambda (\Omega)
   \neq 0$. Suppose $\lambda (\Omega ^c)=\alpha \,(\alpha >0)$.
   By considering another invariant measure $\tilde {\lambda} =(1/\alpha )\lambda
   \vert _{\Omega ^c}$ supported on the complement of $\Omega$, we get $\chi (M)=0$.
   This is a contradiction. So $\lambda (\Omega ^c)=0$
\end{proof}

   Theorem~\ref{thm:4.1} (ii) says for instance that the holonomy group of even dimensional hyperbolic
   manifold can not have a finitely additive invariant probability measure. (But it does have complex
   invariant probability measure.) A much broader class of convex projectively flat manifolds
   should have the same property. And in this case the holonomy group of such manifolds can not
   be amenable since amenability enables one to construct an invariant probability measure
    starting from any probability  measure by averaging process. But in general the holonomy group of
    projectively flat manifold is far from being amenable even when it has an invariant probability measure.
    The case of amenable holonomy group is an interesting special case
    and we can obtain a sharper result as the following Theorem~\ref{thm:4.2} shows.

\begin{theorem}\label{thm:4.2}
 Let $M$ be an even dimensional closed projectively flat manifold with amenable holonomy group.
 Then the followings are equivalent.
\begin{enumerate}
 \item [\rm(i)] The developing map is not onto.
 \item [\rm(ii)] $\chi(M)=0$.
 \item [\rm(iii)] There exists finitely additive invariant probability Borel
  measure $\lambda $
 on $\mathbb R P^n$ such that $\lambda (K) =0$ for any compact subset $K $ of
 the developing image.
 \item [\rm(iv)] For any finitely additive invariant probability Borel measure
 $\lambda $ on
 $\mathbb R P^n$, $\lambda (K) =0$ whenever $K$ is a compact subset of the
 developing image.
 \item [\rm(v)] There exists a countably additive invariant probability Borel measure
 $\lambda $ on $\mathbb RP^n $ such that $\lambda (\Omega )=0$, where $\Omega$ is the developing image of $\tilde M$.
 \item [\rm(vi)] For any countably additive invariant probability Borel measure
 $\mu $, $\mu (\Omega)=0 $.
\end{enumerate}
\end{theorem}
\begin{proof} (i)$\Rightarrow $(ii);
   Let ${\lambda}_0$ be the Dirac measure concentrated at a point $x_0$
    outside the developing image.
   Let $m$ be an invariant mean on $B(H)$.
   Then we can define a measure $\lambda$ on $\mathbb RP^n$ by
    $\lambda(E)=m(f_E)$ for each subset $E \subset \mathbb RP^n$, where $f_E$
    is defined as the equation~\eqref{eq:4.1}.
   Then $\lambda$ is an invariant finitely additive probability measure and
    for each subset $E$ contained in the developing image, $\lambda(E)=0$.
   Therefore $\chi(M)=0$ by the Main Theorem.

(ii)$\Rightarrow$(i) has already been shown in Theorem 4.1 (iii).

    Since the holonomy group is amenable, there exists a finitely additive
     invariant probability Borel measure by averaging and furthermore there
     also exists a countably additive invariant probability Borel measure
     by compactness of $\mathbb RP^n $.

(ii)$\Rightarrow$(iv),(vi);
    Suppose that $\lambda _1$ is a finitely additive invariant probability
     Borel measure on $\mathbb RP^n$ and $\lambda _2$ is a countably additive
     invariant probability Borel measure.
    Let $\mu_1$ and $\mu_2$ be the corresponding induced measure on $M$
     respectively. Then $\mu_1 (M)=\mu_2 (M)=0$ since $\chi (M)=0$.
    By the definition of the induced measure, for each $x$ in the developing
     image, there exist an open neighborhood $U_x$ such that
     $\lambda _1(U_x) =\lambda _2(U_x) =0$.
    Therefore $\lambda _i (K)=0 \, (i=1,2)$ for any compact subset $K$ of
     the developing image and futhermore $\lambda _2(\Omega )=0$ since
     $\lambda _2$ is countably additive.

(iv)$\Rightarrow$(iii) and (vi)$\Rightarrow$(v) are true since the
holonomy
    group is amenable.

(iii) and (v) each imply (ii) by the Main Theorem.
\end{proof}

   Another interesting special case in which the existence of $H$-invariant
    probability measure is guaranteed is where the holonomy group $H$ has
    a fixed point or more generally has a finite orbit.
   In this case, we have the following theorem.

\begin{theorem}\label{thm:4.3}
 Let $M$ be an even dimensional closed projectively flat manifold with holonomy group $H$.
 Suppose $H$ has a finite invariant set $I$.
 Then we have the followings:
\begin{enumerate}
 \item [\rm(i)]  $I\subset D(\tilde M)$ if and only if $\chi(M)> 0.$
 \item [\rm(ii)] $I\cap D(\tilde M)=\varnothing $ if and only if $\chi(M)=0.$
\end{enumerate}
   In particular, if $H$ has a fixed point outside the developing image, then
    the Euler characteristic of $M$ must vanish.
\end{theorem}
\begin{proof}
 Define an invariant probability Borel measure $\lambda $ on $\mathbb RP^n$
  by $\lambda (E) =\sum _{a\in E\cap I} \frac 1n $ for any subset $E$ of
  $\mathbb RP^n$ where $n$ is the cardinal number of $I$.
 Let $\mu $ be the induced measure on $M$.
 Let $I_1 =I \cap D(\tilde M)$ and $I_2 =I\setminus I_1$.
 Suppose that neither $I_1$ nor $I_2$ is empty.
 Let $\lambda_1 =\frac 1{\lambda (I_1)} \lambda |_{I_1}$ and
     $\lambda_2 =\frac 1{\lambda (I_2)} \lambda |_{I_2}$ and $\mu_1$
     and $\mu_2$ are the corresponding induced measures on $M$, respectively.
 Then $\mu_1 (M)>0$ and $\mu_2 (M)=0$. This is a contradiction. Therefore
   either $I_1 =\varnothing $ or $I_2 =\varnothing $.
 If $I_1 =\varnothing $, i.e, $I\cap D(\tilde M)=\varnothing$, then $\chi (M)
 =\mu_2(M)=0$.
 Otherwise, i.e, $I\subset D(\tilde M)$, then $\chi (M)=\mu(M)=\mu_1 (M)>0$.
 Since $\chi (M)\geq 0$, the converses are proved immediately.
\end{proof}

  For the case of an affinely flat manifold $M$, a fixed point can not lie
   in the developing image by the result of Fried, Goldman and Hirsch \cite{3}
   and hence  $\chi (M)$ vanishes if holonomy group $H$ has a fixed point, that is, if $M$ is a radiant affine manifold.
  The vanishing of $\chi (M)$ was observed by Kobayashi \cite{11} using the
   Euler vector field, which gives a non-vanishing $H$-invariant vector
   field on $D(\tilde M)$ and thus a non-vanishing vector field on $M$.

  If $M$ is a closed affinely flat manifold, we can go further to
  show the following Theorem~\ref{thm:4.4} giving an affirmative answer for
  the Chern conjecture when the holonomy group of $M$ has an
  invariant finitely additive probability measure on $\mathbb RP^n$.
  In fact, we do not know whether an affine manifold always have
  such a measure on $\mathbb RP^n$. Anyways, the theorem generalizes
  the earlier result of Hirsch and Thurston for the amenable
  holonomy case \cite{6} and of Kobayashi for the radiant case in a
  unified way. If one can show directly that the holonomy group of a
  complete affine manifold has an invariant probability measure on
  $\mathbb RP^n$, then the theorem would also cover the result of Kostant
  and Sullivan. In fact, if Auslander conjecture is true, that is,
  if the fundamental group of a complete closed affinely flat manifold is
  virtually solvable, then this is an amenable case and
  has an invariant probability measure.

  The holonomy group $H$ of affinely flat manifold of $M$ acts on $\mathbb E^n$
   as affine transformations. Recall that $\mathbb E^n$ is given by
   $x_{n+1} = 1$ in $\mathbb R^{n+1}$.
  The linear parts of these affine transformations are well-defined and
   form a group  called the linear holonomy  group.
  If we projectivize the linear holonomy group, we obtain an action of the
   projectivized linear holonomy group on the projective space, denoted by
   $\mathbb RP^{n-1}$, of the vector space $\mathbb R^n$ associated to the affine
   space $\mathbb E^n$.
  From the Main Theorem, we see immediately that if the projectivized linear
   holonomy group has an invariant probability Borel measure on
    $\mathbb RP^{n-1}$ then $\chi (M)=0$, since such a measure can be regarded
    as a holonomy invariant probability Borel measure on $\mathbb RP^n$ supported
    on $\mathbb RP^{n-1}_{\infty}$ which is disjoint from the affine space $\mathbb E^n$.

\begin{theorem}\label{thm:4.4} Let $M$ be an even dimensional closed affinely
flat manifold with holonomy
 group $H$. If $H$ has an invariant finitely additive probability measure on $\mathbb RP^n$,
 the compactification of $\mathbb E^n$, then $\chi (M)=0$.
\end{theorem}

\begin{proof} Note that there exists a countably additive
$H$-invariant probability measure on $\mathbb RP^n$ by
Propositions \ref{prop:1} and \ref{prop:2} in Appendix. Consider
$M$ as a projectively flat manifold so that $H \subset$
Aff$(n,\mathbb R)\subset$ PGL$(n+1,\mathbb R)$.
 By Furstenberg Theorem \cite[Cor 3.2.2]{16},
 either (i) $\overline H$ is compact or (ii) there
is a proper subspace $V_0$ such that $\lambda \lbrack V_0 \rbrack >0$
and $V_0$ is invariant by a subgroup of $H$ with finite index. An
affinely flat manifold also can be viewed as a $({\mathbb S}^n$,
P$^+$GL$(n+1,\mathbb R))$ manifold, where P$^+$GL$(n+1,\mathbb
R)\cong$ GL$(n+1,\mathbb R)/\mathbb R ^+$.
 Let SL$^{\pm}(n+1,\mathbb R)=\lbrace \,
 A \in$ GL$(n+1,\mathbb R)\,|\, \det A=\pm 1 \, \rbrace$. Then
 P$^+$GL$(n+1,\mathbb R) \cong$ SL$^{\pm}(n+1,\mathbb R)$. Notice that
PGL$(n+1,\mathbb R)\cong$ SL$(n+1,\mathbb R)$ if $n$ is even.
 Let $q\,:\,$ SL$^{\pm}(n+1,\mathbb R) \to$ PGL$(n+1,\mathbb R)$ be the covering homomorphism
 and $p\,:\,{\mathbb S}^n \to \mathbb RP^n$ be the usual covering map.
 Let $D\,:\, \tilde M \to \mathbb RP^n$ be the developing map and
 $\tilde D\,:\, \tilde M \to {\mathbb S}^n$ be its lifting so that $\tilde D \circ p =D$.
 Let $\tilde H \subset$ P$^+$GL$(n+1,\mathbb R)$ be the holonomy group corresponding to $\tilde D$ so that it is the lifting of $H$.

 If $\overline H$ is compact then
 $q^{-1}(\overline H)$ is compact in GL$(n+1,\mathbb R)$ since SL$^{\pm}(n+1,\mathbb R)$
 is closed.
 Therefore there exists $q^{-1}(\overline H)$-invariant inner product on $\mathbb
 R^{n+1}$ and thus we may assume via conjugation that there exists a $\tilde H$-invariant
Riemmanian metric $\phi $ on ${\mathbb S}^n$.
 Since $\tilde D^* \phi $ is deck transformation invariant Riemmanian metric
on $\tilde M$, there is a Riemmanian metric $\psi $ on $M$ such
that the covering map becomes a local isometry. Since $M$ is
compact, $\psi $ is complete and $\tilde D^* \phi$ is also
complete and thus $\tilde D$ becomes a covering map.
 Therefore $\tilde M$ is homeomorphic to ${\mathbb S}^n$ and thus
$M$ is a spherical space form.
 But an affinely flat manifold cannot be a spherical space form.
 So there is a proper subspace $V_0$ such that $\lambda \lbrack V_0 \rbrack >0$
 and $V_0$ is invariant by a subgroup of $H$ with finite index.
 Let $V_0$ be of minimal demension among all linear subspace with
$\lambda \lbrack V_0\rbrack >0$ and $V_0$ is invariant by a subgroup
of $H$ with finite index.
 We may assume $\lbrack V_0\rbrack $ is invariant by $H$.
 If $V_0 \cap \mathbb E^n =\varnothing$ (recall that $\mathbb E^n$ is the affine space
 given by $x_{n+1} =1 $ in $\mathbb R^{n+1}$), then $\chi (M)=0$ by considering an
$H$-invariant probability measure $\lambda' =( 1/{\lambda \lbrack
V_0\rbrack}) \lambda\vert _{\lbrack V_0\rbrack}$.
 Now assume that $V_0 \cap \mathbb E^n \neq \varnothing$.
 If $\dim V_0 =1$, then $V_0 \cap \mathbb E^n$ is a fixed  point  and thus $M$
is radiant.
 Assume $\dim V_0 \geq 2$.
 Consider $\lbrack V_0\rbrack $ and $H' =H\vert_{\lbrack V_0\rbrack}$.
 Again by Furstenberg Theorem, either $\overline {H'}$ is compact or there
exist a proper subspace $W$ of $V_0$ such that $\lambda \lbrack
W\rbrack >0$ and $\lbrack W\rbrack $ is invariant by a subgroup of
$H'$ with finite index.
 But by minimality of $ V_0$, $\overline {H'}$ is necessarily compact
 and thus $q^{-1}(\overline {H'})$ is compact in GL$(m,\mathbb R)$ where $m=\dim V_0 $.
 So there exists a $q^{-1}(\overline {H'})$-invariant inner product on
 $V_0$ and this gives an inner product on $\mathbb R^{n+1}$ such that $q^{-1}
 (H)$ acts  by orthogonal transformation leaving $V_0$ invariant.
 Let $W =V_0 \cap \mathbb E^n$ and  $W_0=V_0 \cap (\mathbb R^n \times \{ 0\})$.
 Observe that $W_0$ is the subspace  of $V_0$ contained in
 $\mathbb R^n \times \{0\}
 $ which is obtained by translating  $W$ along $W_0^\perp$ in $V_0$.
 Then $\lbrack W\rbrack$, $\lbrack W_0\rbrack$ and $\lbrack W_0^\perp \rbrack $
 are all $H$-invariant.
 Furthermore $\lbrack W_0^\perp \rbrack \cap \lbrack W \rbrack
  =\lbrace p\rbrace $ and $p$ is also invariant
and thus $M$ is radiant.
 Since a fixed point can not lie in the developing image in the case of an
affinely flat manifold, $\chi (M)=0$ by considering the Dirac
measure concentrated at $p$.
\end{proof}

\appendix
\section{} In this appendix, we show that for each finitely additive
probability Borel measure on a  compact Hausdorff space there
exists a countably additive probability Borel measure corrsponding
to the measure and furthermore the corresponding countably
additive measure is $G$-invariant if $G$ acts on $X$ and the
finitely additive measure is $G$-invariant.
\begin{proposition}\label{prop:1} Let $\mu_f$ be a finitely additive probability Borel measure on
a compact Hausdorff space $ X$ . Then there exists a countably
additive  probability Borel measure $\mu_c$ on $X$, which
corresponds to the measure $\mu_f$.
\end{proposition}
\begin{proof}Let $B( X, \Sigma)$ be the Banach space consisting
    of all uniform limits of finite linear combination of
    characteristic functions of sets in Borel algebra $\Sigma$. Then
    the dual space of $B( X, \Sigma)$ is isometrically isomorphic to
    the Banach space $ba( X, \Sigma)$ consisting of all bounded
    finitely additive measures on $\Sigma$ (See Theorem IV.5.1 in
    \cite{1}). In this correspondence, a probability measure $\mu_f$ in
    $ba(X, \Sigma)$  corresponds to a positive linear functional
    $\Lambda_{\mu_f} $ on $B(X,\Sigma)$ and
    $\Lambda_{\mu_f}(\chi_{_X})=1$ for the characteristic function
    $\chi_{_X} \in B(X,\Sigma)$. Since $\Sigma$ is the Borel algebra
    on $X$ and $B(X,\Sigma)$ is complete with respect to the supremum
    norm, the Banach space $C(X)$ consisting of all continuous
    functions on {\it compact} space $X$ is a Banach subspace of
    $B(X,\Sigma)$. So the restriction $\Lambda_{\mu_f}|_{C(X)}$ of
    $\Lambda_{\mu_f}$  is a positive linear functional on $C(X)$ with
    $\Lambda_{\mu_f}|_{C(X)}(\chi_{_X})=1$ since $\chi_{_X}\in C(X)$.
    Consequently, we have a countably additive probability Borel measure $\mu_c$ on
    $X$ corresponding to $\Lambda_{\mu_f}|_{C(X)}$ by the Riesz
    Representation Theorem. This completes the proof.
\end{proof}
\begin{remark}
    This correspondence does not imply that $\mu_f(E) = \mu_c(E)$
    for all subset $E$ of $X$ which is contained in the Borel
    algebra. For example, consider a finitely additive translation invariant probability Borel measure $\mu_f$. In fact,
    $\mu_f$ can be regarded as a finitely additive probability
    Borel measure on the closed interval $[- \infty,+ \infty ]$, the
    compactification of ${\mathbb R}^1$, such that
    $\mu_f(\{- \infty \})=\mu_f(\{+ \infty \})=0$. But for the corresponding countably
    additive probability measure $\mu_c$, $\mu_c(\{- \infty ,+ \infty \})=1$. In fact $\mu_f(I)=0$ for any bounded interval
     $I \subset \mathbb R$ and this implies that $\mu_c(\mathbb R) = 0$ using the Monotone Convergence Theorem.
\end{remark}
 \begin{proposition}\label{prop:2}
     Let the group G act on a compact metric space $X$ and $\mu_f$
     be a $G$-invariant finitely additive Borel measure on $X$.
     Then the countably additive probability measure $\mu_c$ which
      corresponds to $\mu_f$ is also $G$-invariant.
 \end{proposition}
\begin{proof}
        $G$-invariance of $\mu_f$ implies that $\mu_f(E) =
    \mu_f(gE)$ for all measurable $E$ and $g \in G$ and it follows that
    $\int f \text{d}\mu_f = \int g\cdot f \text{d}\mu_f$
    for any $f \in C(X)$ where $(g \cdot f)(x)=f(g^{-1}x)$. Since $\mu_f = \mu_c$ on $C(X)$,
    $\int f \text{d}\mu _c = \int g\cdot f \text{d}\mu _c$
    for any $f \in C(X)$, which in turn implies that
\begin{equation*}
\mu_c(E) =\int \chi_E \text{d}\mu_c =\int g^{-1}\cdot \chi_E
    \text{d}\mu_c =\int \chi_{gE}\text{d}\mu_c=\mu_c(gE)
\end{equation*} for all measurable $E$ and $g \in G$ by the
    Monotone Convergence Theorem.
\end{proof}

\end{document}